\documentclass[12pt]{amsart}
\usepackage{amscd}
\usepackage{a4}
\usepackage{amssymb}

\newcommand{\ov}[1]{{\overline{#1}}}

\newtheorem*{thm}{Theorem}
\newtheorem*{prop}{Proposition}
\newtheorem*{lem}{Lemma}
\newtheorem*{cor}{Corollary}
\newtheorem*{defn}{Definition}

\newcommand{\ZZ}{{\mathbb Z}}

\newcommand{\HomX}[3]{{\mathrm{Hom}_{#1}}\:(#2,#3)}

\newcommand{\AnnX}[2]{{\mathrm{Ann}_{#1}}\:(#2)}
\newcommand{\Ker}[1]{{\mathrm{Ker}}\:(#1)}

\renewcommand{\Im}[1]{{\mathrm{Im}}\:(#1)}

\newcommand{\lra}{\longrightarrow}
\newcommand{\dsum}{\displaystyle\sum}
\newcommand{\AH}{A\# H}

\begin{document}
\title{Integrals in Hopf algebras over rings}
\author{Christian Lomp}
\date{\today}
\address{ Departamento de Matem\'{a}tica Pura ,
Faculdade de Ci\^{e}ncias, Universidade do Porto, R.Campo Alegre 687,
4000 Porto, Portugal} \email{clomp@fc.up.pt}
\dedicatory{in memory of my uncle Morris Morgan}
\thanks{I would like to dedicate this work to my uncle Morris Morgan who passed away very suddenly.
Parts of the material presented here are included in the author's doctoral thesis at
the Heinrich-Heine Universit\"at D\"usseldorf. The author would like to express his gratitude to his supervisor
Professor Robert Wisbauer for all his help advice and encouragement.
The author would also like to thank his colleague Lu{\'\i}s Ant{\'o}nio Oliveira for stimulating and helpful discussions about integrals in semigroup rings.
This work is supported by Funda\c{c}{\~a}o para a Ci{\^e}ncia e a Tecnologia (FCT) through the Centro
de Matem{\'a}tica da Universidade do Porto (CMUP). Available as a PDF file from
http://www.fc.up.pt/cmup.}


\begin{abstract}
Integrals in Hopf algebras are an essential tool in studying
finite dimensional Hopf algebras and their action on rings. Over
fields it has been shown by Sweedler that the existence of
integrals in a Hopf algebra is equivalent to the Hopf algebra
being finite dimensional. In this paper we examine how much of
this is true for Hopf algebras over rings. We show that over any
commutative ring $R$ that is not a field there exists a Hopf
algebra $H$ over $R$ containing a non-zero integral but not being
finitely generated as $R$-module. On the contrary we show that
Sweedler's equivalence is still valid for free Hopf algebras or
projective Hopf algebras over integral domains. Analogously for a
left $H$-module algebra $A$ we study the influence of non-zero
left $\AH$-linear maps from $A$ to $\AH$ on $H$ being finitely
generated as $R$-module. Examples and application to separability
are given.
\end{abstract} \maketitle
\section{Introduction}
Alfred Haar introduced a measure $\mu$ on the space of
representable functions $\mathcal{R}(G)$ of a locally compact
group $G$ (see \cite{Haar}). The map $f \mapsto \int fd\mu$ is an
integral $I\in (\mathcal{R}(G))^\ast$. Hochschild exhibited the
Hopf algebra structure of $\mathcal{R}(G)$ and characterised the
$G$-invariance of $I$ as $I$ being a $\mathcal{R}(G)$-colinear map
(see \cite[pp 27-28]{Hochschild}). Finally Sweedler in
\cite{Sweedler69} carried this notion over to arbitrary Hopf
algebras $H$ as $H$-colinear maps $I:H \rightarrow k$ and
characterised them via the augmentation, i.e. $J\star
I(f)=J(1)I(f)$ for all $J\in H^\ast$. Here $J\mapsto J(1)$ is an
augmentation map for $H^\ast$. Sweedler boiled down the definition
of an integral to an element $t$ of an augmented $R$-algebra $A$
such that $at=\epsilon(a)t$ holds for all $a\in A$ (here
$\epsilon: A\rightarrow k$ denotes the augmentation and $R$ a
commutative ring). The existence of such an integral in a Hopf
algebra $H$ is secured if $H$ is finite dimensional. This fact has
been extensively used in studying Hopf algebras and their action
on algebras.
 Sweedler also showed, that the converse is true, i.e. he showed that the existence of a left or right integral implies that the Hopf algebra
 in question had to be finite dimensional. We will discuss how much of this is still valid for Hopf algebras over rings.
 In particular we will show that if over every commutative ring $R$ that is not a field there exists always a Hopf algebra  having a
 non-zero integral but not being finitely generated as $R$-module.

 We show that over a local ring  or an integral domain $R$
 every Hopf algebra $H$ that is projective as $R$-module having a non-zero left or right integral $t$
 must be finitely generated as $R$-module. If $t$ has zero annulator in $R$ then $R$ might be also semiperfect or
 a finite direct product of integral domains. In the third section we characterise separable Hopf algebras.

In the sequel $R$ will denote a commutative ring with unit.
Let $H$ be an augmented $R$-algebra with augmentation (=algebra homomorphism) $\epsilon: H\rightarrow R$.
The base ring $R$ becomes a left and right $H$-module given by $hr := \epsilon(h)r =: rh$ for any $h\in H$ and $r\in R$.
For every left $H$-module $M$ one defines the $R$-submodule of
$\epsilon$-invariants as $$M_{\epsilon}:=\{m\in
M\mid hm = \epsilon(h)m \forall h\in H\}.$$ It is an easy exercise
to check that the map $$\varphi: \HomX{H-}{R}{M} \rightarrow
M_{\epsilon}$$ with $\varphi(f):=f(1)$ is an isomorphism of
$R$-modules and yields that the functors $\HomX{H-}{R}{-}$ and
$(-)_{\epsilon}$ are equivalent. Analogously we can consider $\epsilon$-invariant elements in right $H$-modules and obtain an
equivalence of functors $\HomX{-H}{R}{-}$ and $(-)_{\epsilon}$.

For a Hopf algebra $H$ over $R$ $\epsilon$-invariant elements are called $H$-invariant elements and the submodule of those elements of
a left $H$-module $M$ is denote by $M^H$ rather than $M_{\epsilon}$. $H$-elements of $H$ seen as left (resp. right) $H$-module
are called left (resp. right) integrals and are denoted by $\int_l$ (resp. $\int_r$).
Hence we have  $\HomX{H-}{R}{H}\simeq \int_l$ and
$\HomX{-H}{R}{H}\simeq \int_r$. In the sequel we will make no distinction in
considering integrals as elements in the Hopf algebra or as $H$-linear maps from $R$ to $H$.

We will call a Hopf algebra $H$ over $R$ {\it free, projective, flat,
finite}
resp. {\it non-finite} provided  $H$ is a free, projective, flat,
finitely generated resp. not finitely generated $R$-module.

\section{Integrals as finiteness condition}

\subsection{}
In this section we will study integrals in Hopf algebras as
finiteness conditions. Sweedler proved that Hopf algebras over
fields having a non-zero integral if and only if they are finite
dimensional. This is obviously not anymore true when the ground field is
replaced by a ring. For example take $R=k\times k$ for some field $k$
and let $H_1$ and $H_2$ be $k$-Hopf algebras then $H=H_1\times H_2$ is a
Hopf algebra over $R$ with component-wise operations,
whose submodule of left integrals is equal to
$\int=\int_1 \times \int_2$ where $\int_j$ is the space of integrals of
$H_j$. Hence $\int \neq \{(0,0)\}$ if and only if one of the
Hopf algebras $H_j$ is finite dimensional. An extreme instance of this
situation is by taking  $H_1 = k$ and $H_2 = k[x]$. Then $H=k\times
k[x]$ is a projective, non-finite Hopf algebra over
the semisimple ring $R=k\times k$, that admits a non-zero left integral
$t=(1,0)$.

\subsection{}
We will show that over any ring $R$, that is not a field, there
exists always a non-finite Hopf algebra $H$ that admits a non-zero
left integral.

\begin{prop} For every commutative ring $R$ that is not a field,
there exists a non-finite Hopf algebra over $R$ that contains a non-zero
integral.
\end{prop}

\begin{proof} Let $n$ be a non-zero non-invertible element of $R$.
The polynomial ring $R[X]$ has a Hopf algebra structure over $R$
by letting $X$ be a primitive element, i.e. $\Delta(X)=1\otimes X
+ X\otimes 1$, $\varepsilon(X)=0$ and $S(X)=-X$. The ideal
$I:=<nX>$ is a Hopf ideal of $R[X]$. Let $H:=R[X]/I$ and denote
the elements of $H$ as images $\bar{f(X)}$ of elements $f(X)$ of
$H$ under the canonical projection. Then $H$ is a Hopf algebra
over $R$ with the induced structure maps $\bar{\Delta},
\bar{\varepsilon}$ and $\bar{S}$. Moreover $\bar{n}$ is a non-zero
left and right integral of $H$ since for all $\bar{f(X)} \in H$
with $f(X)=a_0 + Xg(X)$ we have
$$\bar{n}\bar{f(X)} = \bar{na_0} + \bar{nXg(X)} = a_0\bar{n} =
\bar{\varepsilon{f(X)}}\bar{n}.$$ On the other hand $H$ is not
finitely generated as $R$-module since $a$ is not invertible in
$R$.
\end{proof}

Rephrasing above Proposition we can characterise fields in terms
of integrals in Hopf algebras: {\it A commutative ring $R$ is a
field if and only if every Hopf algebra over $R$ that contains a
non-zero left or right integral is finite.}
\subsection{}
If $R$ is decomposable, we can construct a projective non-finite
Hopf algebra over $R$ possessing a non-zero integral.

\begin{cor}\label{ProjectiveCounterExample}
For any commutative ring $R$ with a non-trivial idempotent
there exists a projective non-finite Hopf algebra over $R$ that contains a
non-zero integral.
\end{cor}

\begin{proof} Let $e$ be a non-trivial idempotent. Take $n:=e$ in the proof of the above Proposition.
The ideal $I=<eX>$ is a direct summand of $H$ as $R$-module since
$R[X]=I\oplus Re \oplus R(1-e)[X]$. Thus $H=R[X]/I$ is a
projective $R$-module.
\end{proof}

\subsection{}
Having seen, that the existence of non-zero integrals does not always
imply the finiteness of a Hopf algebra we are going to examine for which
Hopf algebras Sweedler's conclusion is still valid.

Following his proof we can show an analogous
statement for free Hopf algebras. The basic ingredient of
Sweedler's proof is the following Lemma:

\begin{lem} Let $H$ be a projective Hopf algebra over $R$.
If $H$ admits a non-zero left or right ideal that is finitely generated as $R$-module then there exists an
ideal $I$of $R$ such that $IH$ is a finitely generated $R$-module.
\end{lem}

\begin{proof} Let $K$ be a left ideal of $H$ with $K$ finitely generated as $R$-module.
Then $L:=K\leftharpoonup H^*$ is a left $H$-Hopf module and still finitely generated as $R$-module. Note
that $L^{coH} = I1_H$ for some ideal $I$ of $R$. By the fundamental theorem for Hopf modules we get
$L\simeq HL^{coH} = HI$ as Hopf modules. Thus $IH$ is a finitely generated $R$-module.
\end{proof}

Now we deduce easily the following Theorem

\begin{thm}\label{Sweedler_Free} The
following statements are equivalent for a free Hopf algebra $H$.
\begin{enumerate}
\item[(a)] $H$ is finite.
\item[(b)] $H$ contains a non-zero left or right ideal $I$ with $I_R$ finitely generated.
\item[(c)] $H$ contains a non-zero left or right integral.
\end{enumerate}
\end{thm}
\begin{proof} $(a)\Rightarrow (c)$ follows by a result of Pareigis (see \cite{Pareigis71}.\\
$(c)\Rightarrow (b)$ is obvious since $I:=Ht = Rt$ is a non-zero left ideal of $H$ that is cyclic as
$R$-module for every non-zero left integral $t$.\\
$(b)\Rightarrow (a)$ By the Lemma, there exists an ideal $I$ of $R$ such that $IH$ is finitely generated as
$R$-module. Since $H$ is free, there exists an index set $\Lambda$ such that $H\simeq R^{(\Lambda)}$. Thus
$IH\simeq IR^{(\Lambda)} = I^{(\Lambda)}$ being finitely generated implies $|\Lambda|$ being finite.
\end{proof}
\subsection{}
As a corollary we get:

\begin{cor} Let $H$ be a projective Hopf algebra over a semiperfect ring $R$.
Then $H$ contains a non-zero left or right integral
$t$ with $\AnnX{R}{t}=0$ if and only if $H$ is finite.
\end{cor}

\begin{proof} Let $R=R_1\times \cdots \times R_n$ be a finite direct product of local rings. $H$ can
also be decomposed into $H=H_1\times \cdots \times H_n$ where each
$H_i$ is a Hopf algebra over $R_i$. As $H$ is a projective
$R$-module, $H_i$ is a free $R_i$-module. For the unit $e_i$ of
$R_i$ we have $t_i:=e_it \neq 0$ as $\AnnX{R}{t}=0$. Hence $H_i$
contains a non-zero left (or right) integral $t_i$ and must be
finitely generated as $R_i$-module by the above Theorem. Thus $H$
is finitely generated as $R$-module.

On the other hand if $H$ is finite then by a theorem of Pareigis
\cite{Pareigis71} $\int_l$ is a free rank $1$ $R$-module. Thus
there exists a left integral $t$ such that $\AnnX{R}{t}=0$.
\end{proof}

From the example $H=k\times k[x]$ over $R=k\times k$ with integral
$t=(1,0)$, we see that the condition $\AnnX{R}{t}=0$ cannot be
dropped.

\subsection{}
If $R$ is a quasi-Frobenius ring then it is not difficult to show
that a projective Hopf algebra $H$ over $R$ contains a non-zero
left or right integral $t$ with $\AnnX{R}{t}=0$ if and only if $H$
is a cogenerator for all left or right $H$-modules which are
projective as $R$-modules. Moreover if $R$ is semisimple, i.e. a
finite direct product of fields then $H$ is a cogenerator in
$H$-Mod if and only if $H$ is finite.

\subsection{}
Having seen that free Hopf algebras admitting non-zero integrals
must be finite we will examine the situation for projective Hopf
algebras over integral domains.

We will give an obvious generalisation of Sweedler's result.
\begin{thm}\label{SweedlerIntegralDomain} The following statements are equivalent for a projective Hopf algebra $H$
over an integral domain $R$.
\begin{enumerate}
\item[(a)] $H$ is finite.
\item[(b)] $H$ contains a left or right ideal $I$ with $I_R$ finitely
generated.
\item[(c)] $H$ contains a non-zero left or right integral.
\end{enumerate}
\end{thm}

The proof uses Sweedler's theorem and all we need is the following
Lemma:
\begin{lem} Let $M$ be a projective module over an integral domain
$R$ with quotient field $Q$. If $dim_Q (M\otimes_R Q)$ is finite
dimensional, then $M$ is finitely generated as $R$-module.
\end{lem}

\begin{proof} Let $n = dim_Q (M\otimes Q).$ Denote by $R_p$ the localisation of $R$ by an $p\in Spec(R)$.
Since $R$ is an integral domain, we have $R \subseteq R_p
\subseteq Q$ and henceforth :
$$ M \otimes_R Q \simeq (M \otimes_R R_p) \otimes_{R_p} Q =  M_p \otimes_{R_p} Q.$$
As $M$ is projective and $R_p$ is local, $M_p$ is free and there
exists an index set $\Lambda$ such that $M_p \simeq
R_p^{(\Lambda)}$ as $R_p$-module. Thus  $$ M\otimes_R Q \simeq
R_p^{(\Lambda)} \otimes_{R_p} Q \simeq Q^{(\Lambda)}$$ as
$Q$-modules. Since $dim_Q(M\otimes Q)=n$ it follows $|\Lambda|=n$
and $M$ has local constant rank $n$. By a theorem of Vasconcelos
\cite[Proposition 1.3]{Vasconcelos} $M$ is finitely generated as
$R$-module.\end{proof}

\begin{proof}[Proof of the Theorem]
$(a)\Rightarrow (c)$ follows by Pareigis result
\cite{Pareigis71};\\ $(c)\Rightarrow (b)$ is trivial since for any
non-zero left integral $t$ the left ideal $Ht$ is cyclic. Hence we
only have to prove $(b)\Rightarrow (a)$. Let $I$ be a left ideal
of $H$ that is finitely generated as $R$-module. Then $I\otimes_R
Q$ is a non-zero left ideal of $H\otimes_R Q$ and has finite
dimension as $Q$-vector space. Thus by Sweedler's theorem
\ref{Sweedler_Free} $H\otimes_R Q$ is a finite dimensional Hopf
algebra over $Q$. By the Lemma, $H$ is finitely generated as
$R$-module.
\end{proof}

\subsection{} As a corollary we get:
\begin{cor}\label{Cor_prod_domains}
Let $H$ be a projective Hopf algebra over a finite direct product $R$ of integral domains.
If $H$ contains a non-zero left or right integral $t$ with
$\AnnX{R}{t}=0$ then $H$ is finite.
\end{cor}

\begin{proof}
Let $R=R_1\times \cdots \times R_n$ be  a direct sum of integral
domains $R_i$. The Hopf algebra $H$ also has a decomposition
$H=H_1\times \cdots \times H_n$ where each $H_i$ is a projective
Hopf algebra over $R_i$. Since $\AnnX{R}{t}=0, t_i:=e_it\neq 0$
where $e_i$ is the unit of $R_i$. Thus $H_i$ has a non-zero left
or right integral $t_i$ and must be finitely generated as
$R_io$-module by Theorem \ref{SweedlerIntegralDomain} which makes
$H$ finitely generated as $R$-module.
\end{proof}

%

\subsection{}
Since Sweedler boiled down the notion of integrals to the setting
of augmented algebras. We might ask whether there exists integrals
in bialgebras. In particular we will shortly look at integrals in
semigroup rings and we will see that Sweedler's result fails for
bialgebras.

Note that a non-empty subset $I$ of a semigroup $S$ is called a
left ideal if $sI\subseteq I$ for all $s\in S$.


\begin{thm}\label{IntegraleFuerHalbgruppen}
The submodule of left integrals of the bialgebra $R[S]$ is a cyclic left and right ideal of $R[S]$ whose generator
(if non-zero) is of the form $t=\sum_{x\in I} x$ for some (for any) left ideal $I$ of $S$ which is a finite group.
Moreover the $R$-submodule of left integrals is spanned (as $R$-module) by all elements of the form
$\sum_{x\in I} x$ where $I$ is a left ideal of $S$ which is a finite group.
\end{thm}

\begin{proof}
Let $t$ be a left integral of $R[S]$. Write $t=\sum_{s\in K} r_s
s$ where $K=sup(t) = \{s \in S \mid r_s\neq 0\}$
is the support of $t$. Since $st=t$ for all
$s\in S$ we get $$sK = sup(st) = sup(t)=K.$$ Hence $K$
is a finite left ideal of $S$ such that every element of $S$ acts
as a permutation. Choose a finite non-empty subset $I$ of $K$ such
that $sI=I$ for all $s\in S$ and $I$ is minimal with that
property. Since for all $x\in I$ $Ix$ is also a finite non-empty
left ideal of $I$ and $|Ix|\leq |I|$ we get by minimality of $I$
$I=Ix$. Thus any element of the semigroup $I$ is divisible by any
other element of $I$ on both the left and the right. By
 \cite[Theorem V.1.4]{Ljapin} $I$ is a group. Let
$t':=\sum_{x\in I} r_x x$. Then $t'$ is again a left integral of
$R[S]$, because $I$ is a left ideal and a group. In particular
$$\sum_{x\in I} r_x x=t'=yt'=\sum_{x\in I} r_x yx$$
for all $y\in I$. Comparing the coefficients we get $r_{yx}=r_x$
for all $x,y \in I$. In particular if $x$ is the neutral element
$e$ of the group $I$ it follows $r_y = r_e = r$ for all $y\in I$.
Hence $t'=r(\sum_{x\in I} x)$.\\
It is not difficult to see that $\overline{K}:=K\setminus I$ is
again a left ideal of $S$ with the property
$s\overline{K}=\overline{K}$. Repeating our argument we find
another left ideal $I'\subseteq \overline{K}$ of $S$ which is a
group and so forth. Eventually we get left ideals $I_1, I_2,
\ldots, I_k$ of $S$ which are finite
groups such that $K$ is the direct union of those $I_j's$.\\
We may right $$t=\sum_{j=1}^k (\sum_{x\in I_j} r_x x) =
\sum_{j=1}^k r_jt_j$$ where $r_j\in R$ and $t_j = \sum_{x\in
I_j}x$. This shows that every left ideal of $R[S]$ is an
$R$-linear combination of elements of the given form.

Now if $I$ and $J$ are two left ideals of $S$ which are finite groups, then for any $x\in J$
we have that $Ix \subseteq J$ is also a left ideal of $S$. $J$ being a group implies $Ix=J$.
Hence (if exists) we may fix a left ideal $I$ of $S$ which is a finite group, such that every
other left ideal $J$ of $S$ which is a finite group is of the form $Ix$ for some $x\in J$.
If $t$ is a left ideal of $R[S]$ we saw that there are left ideals $I_1,\ldots, I_k$ of $S$ which are finite groups
and $r_1,\ldots, r_k \in R$ such that $t=\sum_{j=1}^k r_jt_j$ where $t_j=\sum_{x\in I_j} x$.
We might choose elements $x_j \in I_j$ such that $I_j=Ix_j$ and
$$t=\sum_{j=1}^k r_jt_j = \sum_{j=1}^k r_j\left(\sum_{y\in I} yx_j\right) = t' \left(\sum_{j=1}^k r_jx_j \right) $$
where $t'=\sum_{y\in I} y$. This shows $t\in t'R[S]$. Since every element in $t'R[S]$ is a left integral,
we proved $\int_l = t'R[S]$.
\end{proof}

Hence the semigroup ring $k[S]$ for any monoid $S$ with zero
element $0\in S$, i.e. $s0=0=0s$ for all $s\in S$, contains a
non-trivial left (and right) integral $0$ independently whether
$S$ is finite or infinite.

\begin{cor} If $R[S]$ contains non-zero left integrals and a non-zero right integrals, then $\int_l=\int_r = Rt$
where $t=\sum_{x\in I} x$ and $I$ is a two-sided ideal of $S$ which is a finite group.
\end{cor}

\begin{proof} If $R[S]$ contains a non-zero left integral, then $S$ contains a left ideal $I$ which is a finite group.
Analogously if $R[S]$ contains a non-zero right integral, $S$ contains a right ideal $J$ of $S$ which is a finite group.
Since $I$ is a left ideal $JI\subseteq I$ and since $I$ is a group $JI=I$. By the analogous argument
we get $J=JI$. Hence every left ideal of $S$ which is a group is an ideal of $S$ and there can exist at most one of those ideals.
Hence there exists exactly one element of the form $t=\sum_{x\in I}x$ where $I$ is a left (or right) ideal of $S$ which is a finite group, i.e.
$\int_l = \int_r = Rt$.
\end{proof}

\subsection{}

Since the only left ideal of a group is the group itself we
conclude that the group ring $R[G]$ contains a non-zero integral
if and only if $G$ is finite. This also holds a bit more general
for cancelative monoids. Recall that a monoid $S$ is called right
cancelative if $ba=ca \Rightarrow b=c$ for all $a,b,c\in S$.
If $S$ is right cancellative, then $R[S]$ contains a non-trivial left
integral if and only if $S$ is a finite group, because if
$I$ is a left ideal of $S$ which is a finite group and $x$ is some element of $I$,
the map
$f: S \rightarrow I$ with $f(s)=sx$ is injective. Hence $|S|\subseteq |I|\subseteq |S|$ implies
$S=I$, i.e. $S$ is a finite group.

%

\subsection{} It is easy to construct finite semigroups that do not
admit left ideals which are groups. In particular finite
completely simple semigroups (see \cite{Howie}) have this
property. An easy way to construct examples is as follows: Fix a
field $k$. Let $I$ and $J$ be non-empty sets and let $S:=I\times
J$. Define a multiplication on $S$ as $(i,j)(i',j'):=(i,j')$ for
all $i,i' \in I$ and $j,j' \in J$. Then $S$ does not admit any
left (resp. right) ideal unless $|I|=1$ (resp. $|J|=1$) (see
\cite{Howie}).
%

Hence for two sets $I$ and $J$ of two elements, the
5-dimensional bialgebra $k[S^1]$ over $k$ does not posses any
non-zero left or right integral (here $S^1$ denotes the semigroup
$S$ with adjoint neutral element).
%
Moreover the previous example shows that the uniqueness of
integrals in Hopf algebras fails for bialgebras since we can
construct $n+1$-dimensional bialgebras whose space of integrals is
$n$-dimensional for any $n\geq 1$. In case $|I|=1$ and
$|J|=\infty$ we can construct an infinite dimensional bialgebra
that posses non-zero left integrals.

\section{Separable Hopf algebras}

In the introduction we saw that there exists a correspondence between integrals and $H$-linear maps from $R$ to $H$.
We will show now another correspondence between integrals and
homomorphisms from $H$ to $H\otimes H$ that is essential for characterising separable Hopf algebras.

\subsection{} Denote by $C_{H\otimes H}(H)$ the set of $H$-centralising elements
of $H\otimes H$, i.e.
$$C_{H\otimes H}(H) := \left\{\sum_i x_i \otimes y_i \in H\otimes H \mid \sum_i hx_i\otimes y_i = \sum_i x_i \otimes y_ih \:\:\:\forall h \in H\right\}.$$
The evaluation homomorphism
$$\Psi: \HomX{H-H}{H}{H\otimes H} \longrightarrow C_{H\otimes H}(H) \text{ with }
\varphi \mapsto \varphi(1_H)$$ is an isomorphism of $R$-modules.
We have the following correspondence between integrals and
$H$-centralising elements:

\begin{lem}\label{LemmaKorrespondenzCasimir} The following homomorphisms of $R$-modules exist:
$$\begin{array}{rcccccc}
i_l:& \int_l &\longrightarrow &C_{H\otimes H}(H),&
t &\mapsto &(1\otimes S)\Delta (t) \\
p_l:& C_{H\otimes H}(H) &\longrightarrow &\int_l, &\sum_i x_i
\otimes y_i &\mapsto & (1\otimes \varepsilon)\left(\sum_i
x_i \otimes y_i\right)\\
\end{array}$$
where $p_l\circ i_l = id_{\int_l}$.
\end{lem}

\begin{proof}
The image of $i_l$ lies in $C_{H\otimes H}(H)$, as for all $t \in
\int_l$ and $h\in H$:
$$ \Delta(t)\otimes h = \sum_{(h)} \Delta(\varepsilon(h_1)t) \otimes h_2 = \sum_{(h)} \Delta(h_1t) \otimes h_2 = \sum_{(h,t)} h_1t_1 \otimes h_2t_2 \otimes h_3.$$
Hence
$$\sum_{(t)} t_1 \otimes S(t_2)h = \sum_{(h,t)} h_1t_1 \otimes S(h_2t_2)h_3 = \sum_{(h,t)} h_1t_1 \otimes S(t_2)S(h_2)h_3 = \sum_{(t)} ht_1 \otimes S(t_2).$$
This shows $i_l(t) \in C_{H\otimes H}(H)$.\\
We will show now that $p_l(u)$ is a left integral for any $u\in
C_{H\otimes H}(H)$. Let $u=\sum_i x_i \otimes y_i \in C_{H\otimes
H}(H)$. For all $h \in H$:
$$ \sum_i hx_i \otimes y_i = \sum_i x_i \otimes y_ih \in H\otimes H.$$
Applying $(1\otimes \varepsilon)$ to this element we get:
$$ h p_l(u) = \sum_i hx_i\varepsilon(y_i) = \sum_i x_i\varepsilon(y_ih) = \varepsilon(h) p_l(u).$$
Thus $p_l(u)\in \int_l$ has been proven.
\end{proof}

Analogous one proves that there are homomorphisms between $\int_r$
and $C_{H\otimes H}(H)$:
$$\begin{array}{rrclrcl}
i_r:& \int_r &\longrightarrow &C_{H\otimes H}(H) &
t &\mapsto &(S\otimes 1)\Delta (t) \\
p_r:& C_{H\otimes H}(H) &\longrightarrow &\int_r &\sum_i x_i
\otimes y_i &\mapsto & (\varepsilon \otimes 1)\left(\sum_i
x_i \otimes y_i\right)\\
\end{array}$$

where $p_r\circ i_r = id_{\int_r}$.

\subsection{}
Recall some definitions: Let $S\subseteq T$ be an extension of rings. $T$ is called a left (resp. right) {\bf semisimple
extension} of $S$ if every short exact sequence of left $T$-modules that splits as left $S$-modules also
splits as left $T$-modules (see \cite[Def. 1]{HirataSugano}).

$T$ is called a {\bf separable extension} of $S$ if the map
$m: T\otimes_S T \longrightarrow T$ with $m(x\otimes y):=xy$ splits as $T$-bilinear map, i.e. if there exists a $T$-bilinear map
$\varphi : T \longrightarrow T \otimes_S T$ with $m\circ \varphi = id_T$ (see \cite[Def. 2]{HirataSugano}). Equivalently
$T$ is a separable extension of $S$ if and only if there exists a separability idempotent $\sum_i x_i \otimes_S y_i \in
T\otimes_S T$ with $\sum_i x_iy_i = 1$ and $\sum_i tx_i \otimes_S y_i = \sum_i x_i \otimes_S y_it$ for all $t\in
T$. In case $S$ is commutative we also call $T$ a separable $S$-algebra if $T$ is a separable extension of $S$.
It is well-known that any separable extension is a semisimple extension.

\subsection{}
The following characterisation of separable Hopf algebras does not require any condition on the Hopf algebra
as module over its base ring.

\begin{thm}\label{SeparabeleHopfalgebren}
The following statements are equivalent for a Hopf algebra $H$ over $R$:
\begin{enumerate}
\item[(a)] $H$ is a separable $R$-algebra.
\item[(b)] $H$ is a left or right semisimple extension of $R$.
\item[(c)] $R$ is a projective left or right $H$-module.
\item[(d)] there exists a left or right integral $t$ in $H$ such that $\varepsilon(t)$ is invertible in $R$.
\end{enumerate}
If $H$ satisfies one of the statements above, then $H$ is separable and
finitely generated over its  centre $Z(H)$
and hence a PI-algebra.
Moreover if additionally $H$ is a projective $R$-module then
$H$ is also finitely generated as $R$-module.
\end{thm}

\begin{proof}
$(a) \Rightarrow (b)$ follows from \cite[28.5]{wisbauer96};\\
$(b) \Rightarrow (c)$ from \cite[20.5]{wisbauer96} it follows that $R$
is an $(H,R)$-projective $H$-module, i.e.  every short exact sequence of left $H$-modules
\[\begin{CD}  0 @>>>  X @>>> Y @>>> R @>>> 0\end{CD}\]
that splits as $R$-modules splits as left $H$-modules.
Hence $\varepsilon$ splits in $H$-Mod, i.e.  $R$ is a projective left (resp. right) $H$-module;\\
$(c) \Rightarrow (d)$ follows from $\HomX{H-}{R}{H}\simeq \int_l$ resp. $\HomX{-H}{R}{H}\simeq \int_r$;\\
$(d) \Rightarrow (a)$ Let $t$ be a left integral with
$\varepsilon(t)$ invertible in $R$. Without loss of generality we
might assume $\varepsilon(t)=1$. By the correspondence
\ref{LemmaKorrespondenzCasimir} $t$ gives rise to the
$H$-centralising element
$$ \omega := i_l(t) = (1\otimes S)\Delta (t) = \sum_{(t)} t_1 \otimes S(t_2) \in H\otimes H.$$
Let $\mu:H\otimes H\rightarrow H$ denote the multiplication map, then
$$ \mu(\omega) = \sum_{(t)} t_1S(t_2) = \varepsilon(t) = 1.$$
Hence $\omega$ is a separability idempotent for $H$ over $R$, i.e. $H$ is a separable extension of $R$.
Analogously one checks for a right integral $t$ with $\varepsilon(t)=1$ that
$$\omega = i_r(t) = (S\otimes 1)\Delta(t) = \sum_{(t)} S(t_1) \otimes t_2$$
is a separability idempotent for $H$ over $R$.

If $H$ is a separable extension of $R$ then $H$ is also a separable extension of its centre $Z(H)$ and hence an
Azumaya algebra (see \cite[28.4]{wisbauer96}).
Central Azumaya algebras are finitely generated and projective and hence PI-algebras
(see \cite[28.1]{wisbauer96}). It follows from \cite[28.5]{wisbauer96} that
projective separable extensions are finitely generated.
\end{proof}

\subsection{}\label{GegebenbeispielBialgebrenNichtSeparabel}
The equivalence $(a)\Leftrightarrow (d)$ has been shown in
\cite[Corollary 5.2]{KadisonStolin00} and \cite{Pareigis73V} for
finite projective Hopf algebras that are Frobenius. As we saw
those conditions are not necessary.

Property $(b)$ of the theorem above is often called the Maschke-property and plays a central role in the study of
semisimple Hopf algebras.

\subsection{}
In contrast to separable Hopf algebras
we will show that semisimple Hopf algebras, i.e. Hopf algebras that are semisimple artinian as rings,
are always finitely generated and projective over their base ring.

\begin{prop}\label{halbeinfacheHopfalgebren} A Hopf algebra $H$ over $R$ is a semisimple ring
if and only if $H$ is separable over $R$ and $R$ is a semisimple ring.
In particular every semisimple Hopf algebra is finitely generated projective as a module over its base ring.
\end{prop}

\begin{proof}
If $H$ is a separable extension of a semisimple ring $R$ then $H$ itself is a semisimple ring as separable extensions are semisimple extensions.\\
Assume that $H$ is a semisimple ring, then  $R$ is a projective left $H$-module and by
 \ref{SeparabeleHopfalgebren} $H$ is separable over $R$. Since every ideal $I$ of $R$ is a direct summand of $R$ as left
 $H$-submodule, $I$ is also a direct summand of $R$ as $R$-module, i.e. $R$ is a semisimple ring.

If $H$ is separable over $R$ and $R$ semisimple then $H$ is projective as $R$-module. By  \cite[28.5]{wisbauer96}
$H$ is finitely generated as $R$-module.
\end{proof}

\subsection{}
The fact that semisimple Hopf algebras over fields are finite dimensional was known
(see \cite[Remark after 2.7]{Sweedler69});
 already the existence of a non-zero integral implies this. In general we see that also any semisimple Hopf algebra over a commutative ring $R$
 has to be finitely generated over $R$ without a priori assumptions on $R$ or on $H$ as $R$-module.

\subsection{}
Separable Hopf algebras that are not semisimple can be easily constructed:

Let $G$ be a finite group of order $n\geq 1$. Set
$R:=\ZZ[\frac{1}{n}]$ the localisation of $\ZZ$ by the set
$\{ 1, n, n^2, \ldots \}$. Then $H:=R[G]$ is a separable Hopf algebra over $R$ by
\ref{SeparabeleHopfalgebren}, since
$n=|G|=\varepsilon(t)$ (with  $t=\sum_{g \in G} g$) is invertible
in $R$. But $H$ is not semisimple as $R$ is not semisimple.

\subsection{}
Let $H$ be a Hopf algebra over $R$. A pure $R$-submodule $K$ of $H$ is called a
Hopf subalgebra if $K$ is a subalgebra of $H$ with $\Delta(K)\subseteq
K\otimes K$ and $S(K)\subseteq K$. A Hopf subalgebra $K$ of $H$ is called normal if
$K$ is closed under left and right adjunction, i.e. if for all $h\in H$ and $k\in K$:
$$\dsum_{(h)} h_1kS(h_2) \in K
\:\mbox{ und }\: \dsum_{(h)} S(h_1)kh_2 \in K.$$
For a normal Hopf subalgebra $K$ the quotient algebra $\ov{H}=H/K^+H$ with
$K^+:=K\cap \Ker{\varepsilon}$ becomes a Hopf algebra over $R$. Moreover we get also $K^+H = HK^+$.
The following theorem extends  \cite[Theorem
3.10]{Linchenko} from Hopf algebras over fields to Hopf algebras over rings:

\begin{prop}\label{SeparabelDurchUnterUndQoutienten}
Let $H$ be a Hopf algebra over $R$ and let $K$ be a normal Hopf subalgebra of $H$
with quotient $\ov{H}:=H/K^+H$ where $K^+=K\cap
\Ker{\varepsilon}$. Suppose that $K$ and $\ov{H}$ are separable $R$-algebras, then $H$ is also a separable $R$-algebra.
\end{prop}

\begin{proof}
Let $t$ be a left integral in $K$ with  $\varepsilon(t)=1$ and let
$\ov{s}$ be a left integral in $\ov{H}$ with $\ov{\varepsilon}(\ov{s})=\varepsilon(s)=1$.
Then we also have $\varepsilon(st)=1$. We show now that $st$ is a left integral in $H$.
Let $h\in H$. From $\ov{h}\ov{s} = \varepsilon(h)\ov{s}$ in $\ov{H}$ it follows that there exists  $x\in K^+H$ with
$$hs-\varepsilon(h)s = x \in K^+H=HK^+.$$ Since $K^+t=0$ we get
$hst -\varepsilon(h)st = xt = 0$. By
\ref{SeparabeleHopfalgebren} $H$ is separable over $R$.
\end{proof}

\subsection{}

\begin{lem}\label{SeparabeImpliziertUnterUndQuotienten}
Let $H$ be a Hopf algebra over $R$ and $K$ a normal Hopf subalgebra.
with quotient $\ov{H}$. If $H$ is a separable $R$-algebra then so is $\ov{H}$
and if moreover $H$ is free as left $K$-module, then $K$ is also a separable $R$-algebra.
\end{lem}

\begin{proof}
We have $\ov{H} = H/K^+H \simeq H \otimes_K R$. Let  $f$ be the ring homomorphism
$$f:=1\otimes \varepsilon: \:H\simeq H\otimes_K K \lra H \otimes_K R.$$
Then $f$ is  surjective. Since $H$ is separable over $R$,  $H$ is separable over $K$ by
 \cite[Prop 2.5(1)]{HirataSugano} and by
\cite[Prop. 2.4]{HirataSugano} $f(H)$ is separable over $f(K)$.
From $f(H) \simeq \ov{H}$ and $f(K)\simeq R$ it follows that  $\ov{H}$ is a separable $R$-algebra.\\
To show that $K$ is separable over $R$ if $H$ is free as left $K$-module one can proceed as in
 \cite[2.2.2(2)]{Montgomery} (the proof there does not require that $R$ is a field or any other properties of $H$ as $R$-module).
\end{proof}

\subsection{}

\begin{prop}\label{SeparableUnterHopfalgebren}
Let $H$ be a free Hopf algebra over a noetherian local ring $R$ and let $K$ be a normal Hopf subalgebra of $H$.
Then $H$ is a separable $R$-algebra if and only if $K$ and $\ov{H}$ are separable $R$-algebras.
\end{prop}

\begin{proof}
If $K$ and $\ov{H}$ are separable over $R$ then $H$ is also separable over $R$ by
\ref{SeparabelDurchUnterUndQoutienten}. On the other hand if $H$ is separable over $R$ it must be finitely generated as
$R$-modules since $H$ is free. As $R$ is noetherian $K$ is finitely generated as $R$-module. By
\cite[Lemma 5.2]{KadisonStolin99} $H$ is a free left $K$-module. Hence  $K$ and $\ov{H}$ are separable $R$-algebras by
\ref{SeparabeImpliziertUnterUndQuotienten}.
\end{proof}



\section{Homomorphisms from a module algebra to its smash product}

As we saw in the last section, the existence of a non-zero
integral can be seen as a finiteness condition for the Hopf
algebra. Integrals and $H$-linear maps from $R$ to $H$ are in
bijection to each other since every left integral $t$ defines a
left $H$-linear map $r\mapsto rt$ between $R$ and $H$ and since
for any left $H$-linear map $f:R\rightarrow H$ the element $(1)f$
is a left integral. We will see that more generally the existence
of a  non-trivial homomorphism between an $H$-module algebra $A$
and its smash product $\AH$ can be seen as a finiteness condition
on $H$.

Let $H$ be a bialgebra over $R$ and let $A$ be an $R$-algebra that
is also a left $H$-module algebra. If the multiplication
$\mu:A\otimes A\rightarrow A$ and the unit $\eta: R\rightarrow A$
are $H$-linear maps, then $A$ is a called a left $H$-module
algebra. The smash product of $A$ and $H$ is the $R$-algebra $\AH$
whose underlying $R$-module is $A\otimes H$ and whose
multiplication is given by
$$ (a\# h)(b\# g) := \sum_{(h)} a(h_1\cdot b)\# h_2$$
for all $a\#h, b\#g \in \AH$ where $\cdot$ denotes the left
$H$-action on $A$.

$\AH$ is an $R$-algebra with subring $A\simeq \{ a\# 1 \mid a\in
A\}\subseteq \AH$. Note that $H$ is in general not a subring of
$\AH$. Necessary conditions for this to happen are that $H$ is a
flat $R$-module and $A$ is faithful as $R$-module.

One checks that
$$\HomX{\AH}{A}{\AH} \simeq (\AH)^H = r.ann_{\AH}(\Ker{\alpha})$$
where the first isomorphism is given by $f\mapsto (1)f$.

\subsection{} Let us begin with an example. Let $A$ be an
$R$-algebra and $S$ a monoid whose elements act as endomorphisms
on $A$. Denote the action of an element $g\in S$ on an element
$a\in A$ by $a^g$. Then $A$ is a $R[S]$-module algebra where we
consider the semigroup ring $R[S]$ as an bialgebra. The smash
product $A\# R[S]$ is called the skew-semigroup ring of $A$ and
$S$ and is denoted by $A*S$. We denote the $S$-invariant elements
of $A*S$ by $(A*S)^S$. Note that $\HomX{A*S}{A}{A*S} \simeq
(A*S)^S$.

In this situation we have the following analogue to
\ref{IntegraleFuerHalbgruppen}:

\begin{thm}\label{IntegraleFuerSchiefgruppenringe}
The submodule of $S$-invariant elements $(A*S)^S$ of $A*S$ is a cyclic right ideal of $A*S$ and
generated (if non-zero) by some element $1*t$ where $t=\sum_{x\in I} x$ and $I$ is a left ideal of $S$ which is a finite group.
\end{thm}

\begin{proof}
The proof goes analogously to the proof of
\ref{IntegraleFuerHalbgruppen}. Let $\gamma = \sum_{x\in S} a_x *
x \in (A*S)^S$ for some $a_x \in A$. Then $I=\{x\in S \mid a_x\neq
0\}$ is a finite left ideal of $S$ with $sI=I$ for all $s\in S$.
As in \ref{IntegraleFuerHalbgruppen} one decomposes $I$ into a
disjoint union of left ideals $I_k$ which are finite groups. For all $y\in S$ we have
$$ \sum_{k=1}^n \sum_{x\in I_k} a_x * x = \gamma = (1*y)\gamma = \sum_{k=1}^n \sum_{x\in I_k} a_x^y * yx.$$
As the $I_k$ are disjoint we must have $$\sum_{x\in I_k} a_x * x =
\sum_{x\in I_k} a_x^y * yx$$ for all $k$. In particular if $y \in
I_k$ for some $k$ we have $a_{xy}=a_x^y$ for all $x\in I_k$.
If $x=e$ is the neutral element of $I_k$ we get $a_y = a_e^y$.
Set $b_k:=a_e$. Hence $$\gamma = \sum_{k=1}^n \sum_{x\in I_k} a_x * x
= \sum_{k=1}^n \sum_{x\in I_k} b_k^x * x = \sum_{k=1}^n \sum_{x\in I_k} (1*x)(b_k*1) $$
Fix any left ideal $I$  of $S$ which is a finite group. As seen in the proof of \ref{IntegraleFuerHalbgruppen}
we may choose elements $y_k \in I_k$ such that $I_k = Iy_k$.
Thus $$ \gamma=\sum_{k=1}^n \sum_{x\in I_k} (1*x)(b_k*1)
=\sum_{k=1}^n \sum_{z\in I} (1*zy_k)(b_k*1)
= \left( \sum_{z\in I} 1*z \right)\left(\sum_{k=1}^n (1*y_k)(b_k*1)\right).$$
For $\gamma' :=\sum_{k=1}^n (1*y_k)(b_k*1)$ and $t:=\sum_{z\in I}z$ we have
$\gamma = (1*t)\gamma'$. Note that $t$ does not depend on $\gamma$, hence $(A*S)^S \subseteq  (1*t)R[S]$.
The converse is clear.
\end{proof}

We see that $\HomX{A*S}{A}{A*S} \neq 0$ if and only if $S$
contains a left ideal which is a finite group.
In case $S$ is a left cancellative monoid we see that
$\HomX{A*S}{A}{A*S}\neq 0$ if and only if $S$ is a finite group.

It is known that if a group $G$ acts on a ring $A$ as
automorphisms such that $A$ is a projective module over the skew
group ring $A*G$ then  $G$ has to be finite. We see that already
the existence of a non-zero $A*G$-homomorphism $A\rightarrow A*G$
makes the group finite.

\subsection{}
We will see that under some additional conditions non-zero
$\AH$-linear maps from $A$ to $\AH$ give rise to the existence of
non-zero left integrals in $H$ (which in turn might imply that $H$
is finitely generated as $R$-module as we had seen before).

First we need the following Lemma. For any $A$-module $M$ denote
by $\langle H \otimes M\rangle$ the  left $H$-module whose
$H$-action is given by $h(g\otimes m) = hg \otimes m$ for all $h,
g\in H$ and $m\in M$.

\begin{lem}\label{ProjektiveUndInvarianten} Let $H$ be a Hopf algebra over $R$ and let $P$ be a left  $H$-module that
is projective as $R$-module. Then $\langle H\otimes P\rangle^H =
\int_l \otimes P$.
\end{lem}

\begin{proof}
By hypothesis there exists an index set $\Lambda$ and a split
epimorphism  $\pi: R^{(\Lambda)} \lra P$. The map
$$1\otimes \pi: \langle H \otimes R^{(\Lambda)}\rangle \lra
\langle H \otimes P\rangle$$ is a splitting left $H$-module
homomorphism and induces a splitting  $R$-linear map
$$ p: \HomX{H}{R}{ \langle H\otimes R^{(\Lambda)} \rangle} \lra \HomX{H}{R}{ \langle H\otimes P \rangle}.$$
Consider the following commuting diagram:
\[
\begin{CD}
{\HomX{H}{R}{\langle H\otimes P\rangle}} @<<p< {\HomX{H}{R}{\langle H\otimes R^{(\Lambda)}\rangle}} @>{\simeq}>> {\HomX{H}{R}{H}^{(\Lambda)}}\\
@V{\varphi_1}VV @V{\varphi_2}VV @V{\varphi_3}VV\\
{\langle H\otimes P \rangle} @<{1\otimes \pi}<< {\langle H\otimes R^{(\Lambda)} \rangle} @>{\simeq}>> {H^{(\Lambda)}}
\end{CD}
\]

where $\varphi_1, \varphi_2$ and $\varphi_3$ are the evaluation
maps that evaluate the homomorphisms at $1$. The images of
$\varphi_i$ are precisely the submodules of $H$-invariant elements. We get:
$$\Im{\varphi_3} = \langle H^{(\Lambda)}\rangle^H = \int_l^{(\Lambda)} $$
and hence
$$\langle H\otimes R^{(\Lambda)}\rangle^H = \int_l \otimes R^{(\Lambda)}.$$
Eventually we have
$$\langle H\otimes P\rangle^H = \Im{p\varphi_1} = \Im{\varphi_2(1\otimes \pi)}
= \left(\int_l \otimes R^{(\Lambda)}\right)(1\otimes \pi) = \int_l \otimes P.$$
\end{proof}


\subsection{}
The next Lemma is known for Hopf algebra actions over fields and
will be important for the rest of this section. Let $A$ be a left
$H$-module algebra and let $\alpha:\AH \longrightarrow A$ be the
$\AH$-linear projection given by $a\#h \longmapsto
a\varepsilon(h)$.

\begin{lem}\label{HomAAH}
Let $H$ be a Hopf algebra over $R$ with bijective antipode. Then
$$\left(\AH\right)^H = \left(1\#\int_l\right)\left(A\#1\right)$$ for every left
$H$-module algebra $A$ that is projective as $R$-module.
\end{lem}

\begin{proof}
By Lemma \ref{ProjektiveUndInvarianten} $\langle H\otimes
A\rangle^H = \int_l \otimes A.$ Moreover it is easy to check that
the map
$$ \phi: \langle H\otimes A\rangle \lra \AH \:\:\mbox{ with }
h\otimes a \mapsto (1\# h)(a\#1) $$ is an isomorphism of left
$H$-modules with inverse $$\phi^{-1}(a\# h):= \sum_{(h)} h_2
\otimes S^{-1}(h_1)\cdot a.$$ Hence
$$(\AH)^H = \phi\left(\int_l \otimes A\right) = \left(1\#
\int_l\right)(A\# 1).$$
\end{proof}

\subsection{}
Now we are ready to state the main result of this section relating
non-zero homomorphisms  between a module algebra and its smash
product with non-zero integrals of the Hopf algebra.

\begin{thm}\label{EinbettungVonANachAH} Let $H$ be a Hopf algebra over $R$ with bijective antipode and let $A$
be a left $H$-module algebra that is projective as $R$-module. The
following statements are equivalent:
\begin{enumerate}
    \item[(a)] $\HomX{\AH}{A}{\AH}\neq 0$
    \item[(b)] $\exists \varphi\in\HomX{H-}{R}{H}$ such that $\varphi\neq 0$ and $\Ker{\varphi}A\neq A$.
    \item[(c)] $\exists t\in \int_l$ such that the map $A\lra \AH$ with $a\longmapsto a\# t$ is non-trivial.
\end{enumerate}
If $\int_l$ is cyclic then  $A$ is isomorphic to a left ideal of
$\AH$ if and only if there exists a non-zero left integral $t$ in
$H$ such that $\AnnX{R}{t}\subseteq \AnnX{R}{A}$.
\end{thm}

\begin{proof}
By \ref{HomAAH} $(\AH)^H = (1\#\int_l)(A\#1)$.\\
$(a)\Rightarrow (c)$ Let $0\neq \phi\in\HomX{\AH}{A}{\AH}$. Then
$(1_A)\phi \in (\AH)^H= (1\#\int_l)(A\#1)$. Hence for all $a\in
A$: $(a)\phi = \sum_{i=1}^k (a\#t_i)(a_i\#1)$ for some non-zero
left integrals $t_i\in \int_l$ and elements $a_i\in A$. Thus there
exists at least one left integral $t$ such that the map $a \mapsto a\#t$ is non-trivial.\\
$(c)\Rightarrow (b)$ Assume there exists a non-zero left integral
$t$ such that the map $R_t:A\rightarrow \AH$ with $(a)R_t = a\# t$
is non-trivial. Let $\varphi:R\longrightarrow H$ be the $H$-linear
map $\varphi(r)=rt$. Obviously $\Ker{\varphi}=\AnnX{R}{t}$ and
$\Ker{\varphi}A\subseteq \Ker{R_t}$. Since $t\neq 0$, $\varphi\neq
0$ and since $\Ker{R_t}\neq A$, $\Ker{\varphi}A\neq A$.\\
$(b)\Rightarrow (a)$ Let $\varphi \in \HomX{H}{R}{H}$ be such that
$\Ker{\varphi}A\neq A$. Set $$\phi:=1_A\otimes \varphi: A\simeq
A\otimes R \longrightarrow \AH.$$ $\phi$ is left $\AH$-linear. As
$A$ is a flat $R$-module and $\Ker{\varphi}A\neq A$,  $\phi \neq
0$.

Assume that $\int_l$ is cyclic. If there exists an embedding
$\phi:A \rightarrow \AH$ as $\AH$-modules, then $(1_A)\phi =
(1\#t)(a\# 1)$. Let $r\in \AnnX{R}{t}$ then $$(r1_A)\phi = r(1\#
t)(a\#1) = 0$$ implies that $r1_A=0$, i.e. $r\in
\AnnX{R}{1_A}=\AnnX{R}{A}$. Thus $\AnnX{R}{t}\subseteq
\AnnX{R}{A}$.

On the other hand let $t$ be a left integral such that
$\AnnX{R}{t}\subseteq \AnnX{R}{A}$. Consider the map $R_t:
R\rightarrow H$ with $R_t(r)=rt$. Since $A$ is flat we get also a
map $\phi:A\rightarrow \AH$ with $(a)\phi = a\# t$ whose kernel is
$\Ker{R_t}A$. Since $\Ker{R_t}=\AnnX{R}{t}\subseteq \AnnX{R}{A}$
the kernel is zero and $\phi$ is injective. It is easy to verify
that $\phi$ is left $\AH$-linear.
\end{proof}

\subsection{}

In general the existence of a non-zero left integral in $H$ is not
sufficient for the existence of a non-zero left $\AH$-linear map
from $A$ to $\AH$. We will construct an example:

Let $e$ be a non-trivial idempotent in $R$ and define
$H:=R[X]/\langle eX\rangle$ as in \ref{ProjectiveCounterExample}. Then $H$ is a
projective Hopf algebra over $R$ with bijective antipode and
non-zero integral. Set $A:=R(1-e)$. Then $A$ is an $R$-algebra and
a left $H$-module algebra with trivial $H$-action given by
$\epsilon$. The smash product is isomorphic to
$$\AH = R(1-e)\otimes (R[X]/\langle eX\rangle) \simeq
R(1-e)[X] = A[X]$$ and we have for the group of homomorphisms:
$$\HomX{\AH}{A}{\AH} \simeq \HomX{A[X]}{A}{A[X]} = 0.$$

Thus if $k$ is a field and $R=k\times k$. Then there exists a
Hopf algebra $H$ with bijective antipode and non-zero integrals
over the semisimple ring $R$ and a left $H$-module algebra $A$
with $\HomX{\AH}{A}{\AH}=0$.

\subsection{}
Combining the last theorem with results from the second section
enables us to state  characterisation of when a Hopf algebra is finitely generated over its base ring
in terms of homomorphisms from a module algebra to its  smash product.

\begin{thm} Let $R$ be an integral domain or a local ring.
The following statements are equivalent for a projective Hopf algebra $H$ over $R$.
\begin{enumerate}
\item[(a)] $H$ is finite.
\item[(b)] Every left $H$-module algebra
$A$ which is projective as $R$-module is isomorphic to a left
ideal of $\AH$.
\item[(c)] $H$ has a bijective antipode and there
exists some left $H$-module algebra $A$ which is projective as
$R$-module such that $\HomX{\AH}{A}{\AH}\neq 0$.
\end{enumerate}
\end{thm}

\begin{proof}
$(a)\Rightarrow (b)$ Let $H$ be a Hopf algebra that is finitely
generated and projective as $R$-module. By \cite{Pareigis71} $H$
has a bijective antipode and contains a non-zero left integral
$t$.\\
If $R$ is local then $R$ has trivial Picard group and
$\int_l\simeq R$ by \cite{Pareigis71}. Applying
\ref{EinbettungVonANachAH} we get that $A$ is isomorphic to a left
ideal of $\AH$ for any left $H$-module algebra $A$ that is
projective as $R$-module.\\
If $R$ is an integral domain and $t$ a non-zero left integral,
then $\AnnX{R}{t}=0$ as $H$ is a torsionfree $R$-module. Thus the
map $\varphi:R \rightarrow H$ with $\varphi(r)=rt$ is an injective
$H$-linear map. Take any left $H$-module algebra $A$ that is
projective as $R$-module. Since $A$ is flat as $R$-module also
$\phi: A \rightarrow \AH$ with $\phi(a)=a\# t$ is injective. One
checks easily that $\phi$ is $\AH$-linear. Hence $A$ is isomorphic
to a left ideal of $\AH$.\\
$(b)\Rightarrow (c)$ Since $A=R$ is a left $H$-module algebra with
trivial $H$-action we get that $\HomX{R\# H}{R}{R\# H} \simeq
\HomX{H}{R}{H} \neq 0$. Thus $H$ contains a non-zero left integral
and must be finitely generated by \ref{Sweedler_Free} resp. \ref{SweedlerIntegralDomain}. Hence
$H$ has a bijective antipode.\\
$(c)\Rightarrow (a)$ Assume that $H$ has a bijective antipode and
that there exists a left $H$-module algebra $A$ that is projective
as $R$-module such that $\HomX{\AH}{A}{\AH}\neq 0$. Then by
\ref{HomAAH}  $H$ contains a non-zero left integral. By
\ref{SweedlerIntegralDomain} resp. \ref{Sweedler_Free} $H$ has to
be finitely generated as $R$-module.
\end{proof}

\subsection{}
Let $R$ be an integral domain or a local ring and $H$ a Hopf
algebra over $R$ with bijective antipode. We have in particular
the following dichotomise:\\
\begin{tabular}{cl}
{\bf either}& every left $H$-module algebra $A$ that is projective
as $R$-module \\& is isomorphic to a left ideal of $\AH$\\
{\bf or}& $\HomX{\AH}{A}{\AH}=0$ for every left $H$-module algebra
$A$\\& that is projective as $R$-module.
\end{tabular}

Thus if $H$ is a Hopf algebra with bijective antipode over a field $k$ and $A$ is
a left $H$-module algebra that is projective as $\AH$-module then $H$ must be finite dimensional.

\subsection{}
In the last part of this chapter we will relate integrals to the
projectivity of a module algebra as module over its smash-product.

\begin{defn}
Let $H$ be Hopf algebra over $R$. An element $a\in A$ of a left
$H$-module algebra $A$ is called a {\bf trace one} element if
there exists a left integral $t\in \int_l$ such that $t\cdot a =
1$.
\end{defn}

Obviously the existence of a trace one element in a $H$-module
algebra requires the existence of a left integral and as we have
seen might imply that $H$ is finitely generated  as $R$-module.

Since $t\cdot A \subseteq A^H$ always holds, one easily shows that
$A$ has a trace one element if and only if $t\cdot A = A^H$ holds.

The relation of the existence of trace one elements in a module
algebra $A$ and the projectivity of $A$ as left $\AH$-module is
given in the following proposition.

\begin{prop}\label{Spur1Element}
Let $H$ be a Hopf algebra over $R$ and let $A$ be a left
$H$-module algebra that has a trace one element, then $A$ is a
projective left $\AH$-module.
\end{prop}

\begin{proof}
Let $t\in \int_l$ and $a\in A$ with $t\cdot a = 1$. The map
$\beta: A \lra \AH$ with $(x)\beta:= (x\#t)(a\# 1)$ is left
$\AH$-linear, since for all $x\in A$ and $h\in H$:
\begin{eqnarray*}
h (x)\beta = (1\# h)(x\# t)(a\# 1)
&=& \left( \dsum_{(h)} (h_1\cdot x)\# (h_2t) \right)(a\# 1)\\
&=& \left( \dsum_{(h)} (h_1\cdot x)\# \varepsilon(h_2)t \right)(a\# 1)\\
&=& ((h\cdot x)\# t) (a\# 1) = (h\cdot a)\beta.
\end{eqnarray*}
Moreover  $\beta$  lets $\alpha$ split as
$$ (x)\beta\alpha = [(x\# t)(a\# 1)]\alpha = x(t\cdot a) = x,$$
for all $x\in A$. Thus $A$ is a projective left $\AH$-module.
\end{proof}

\begin{prop}\label{Spur1Element2} Let $H$ be a Hopf algebra over $R$ with bijective antipode
such that $\int_l$ is a cyclic right ideal of $H$ and let $A$ be a
left $H$-module algebra that is projective as $R$-module. If $A$
is projective as left $\AH$-module then $A$ has a trace one
element.
\end{prop}

\begin{proof}
Assume that $A$ is a projective left $\AH$-module. Then $\alpha$
splits, i.e. there exists a $\AH$-linear map $\beta: A \lra \AH$
such that $(x)\beta\alpha = x$ for all $x \in A$. Since the
isomorphism
$$\HomX{\AH}{A}{\AH} \simeq (\AH)^H=(1\# \int_l)(A\# 1$$ (see
\ref{HomAAH}) maps $\beta$ to $(1)\beta$ we might assume $(1)\beta
= (1\# t_i)(a_i\# 1)$. As $\int_l$ is a cyclic right ideal of $H$
we might write $(1)\beta$ as $(1\# t)(a\# 1)$. Now
$$1= (1)\beta\alpha = ((1\# t)(a\# 1))\alpha = t\cdot a$$
shows that $A$ has a trace one element.
\end{proof}

In particular for a ring $R$ with trivial Picard group and a Hopf
algebra $H$ that is finitely generated projective as $R$-module we
have for all left $H$-module algebras $A$ which are projective as
$R$-module: $A$ has a trace one element if and only if $A$ is
projective as left $\AH$-module.

\section{Separable Smash products}

Let $S\subseteq T$ be an extension of rings. A short exact sequence in $T$-Mod  is called
 $(T,S)$-exact if it splits in $S$-Mod. A left $T$-module $M$ is called $(T,S)$-semisimple if every $(T,S)$-exact sequence
 in $\sigma[_TM]$ splits. (Here $\sigma[_TM]$ denotes the full subcategory of $T$-Mod whose objects are quotients of direct sums of $M$).

If $T$ is a $(T,S)$-semisimple left $T$-module if and only if $T$ is a left semisimple extension of $S$.


\subsection{}
\begin{thm}\label{SeparabeleSmashProdukte}
Let $H$ be a Hopf algebra over $R$ and $A$ a left $H$-module
algebra. Suppose that there exists a right integral $t$ in $H$ and
a central element $z$ in $AS$ such that $S(t)\cdot z=1$, them $A$
is a $(\AH, A)$-semisimple projective left  $\AH$-module.
\end{thm}

\begin{proof}
Set $\omega:=\sum_{(t)} (1\# S(t_1))(z\# t_2) \in \AH$. For all $a\in A$ we have:
$$\omega\cdot a = \sum_{(t)} S(t_1)\cdot (z (t_2 \cdot a)) = \sum_{(t)} (S(t_2)\cdot z)(S(t_3)t_2 \cdot a) = (S(t)\cdot z)a=a.$$
Let $M\in \sigma[_{\AH}A]$, then there exists an index set $\Lambda$ and a $\AH$-submodule $I\subseteq A^{(\Lambda)}$, such that
$M$ is isomorphic to a $\AH$-submodule of $A^{(\Lambda)}/I$. Without loss of generality we identify
$M$ with a submodule of  $A^{(\Lambda)}/I$. Let $m\in M$, then there are elements
 $a_\lambda \in A$ with $m = (a_\lambda)_\Lambda +
I$. Hence
$$\omega \cdot m = (\omega\cdot a_\lambda)_\Lambda + I = (a_\lambda)_\Lambda + I = m \:\:\: (\dag)$$

Consider the $\AH$-bimodule $\AH {\otimes_A} \AH$ and its element
$$\Omega:= \sum_{(t)} (1\# S(t_1)) \otimes (z\# t_2).$$ We will
show that $\Omega$ is $\AH$-centralising: Let $h\in H$ then
$$(1\# h)\Omega = \sum_{(t)} (1\# hS(t_1)) \otimes (z\# t_2) = \sum_{(t)} (1\# S(t_1)) \otimes (z\# t_2h)
= \Omega (1\# h)$$ as $\sum_{(t)} S(t_1)\otimes t_2 = i_r(t) \in
C_{H\otimes H}(H)$ by \ref{LemmaKorrespondenzCasimir}.

Let $a\in A$ then
\begin{eqnarray*}
\Omega (a\# 1)
&=&\sum_{(t)} (1\# S(t_1))\otimes (z\# t_2)(a\#1)\\
&=& \sum_{(t)} (1\# S(t_1))\otimes (z (t_2\cdot a)\# t_3)\\
&=& \sum_{(t)} (1\# S(t_1))(t_2\cdot a)\# 1)\otimes (z \# t_3)\\
&=& \sum_{(t)} (S(t_2)t_3 \cdot a \# S(t_1))\otimes (z \# t_4)\\
&=& \sum_{(t)} (a \# S(t_1))\otimes (z \# t_2)\\
&=& (1\# a) \Omega.
\end{eqnarray*}
Thus $\Omega$ is a $\AH$-centralising  element of $\AH \otimes_A
\AH$.

Let  $f\in \HomX{A}{M}{N}$ for two left $\AH$-modules $M,N$.
Since $\Omega$ is $\AH$-centralising we get that
$${\tilde f}: [m \mapsto \sum_{(t)} (1\# S(t_1)) \cdot f((z\# t_2) \cdot m)]$$
is a left $\AH$-linear map.
Hence if $M,N \in \sigma[_{\AH}A]$ and if $N$ is direct summand of $M$ as left $\AH$-module
with $A$-linear projection.
$\pi: M\longrightarrow N$, then $\tilde{\pi}$ is a $\AH$-linear projection, since for all $n\in N$
holds $\pi(n)=n$. It follows from $(\dag)$:
$$\tilde{\pi}(n) = \sum_{(t)} (1\# S(t_1)) \cdot \pi( (z\# t_2) \cdot n)
= \sum_{(t)} (1\# S(t_1)) \cdot ( (z\# t_2) \cdot n) = \omega \cdot n = n.$$
Hence every short exact sequence in $\sigma[_{\AH}A]$ that splits in $A$-Mod also splits in $\AH$-Mod.
\end{proof}


\subsection{}
In case $H$ has a bijective antipode then the hypothesis of the theorem above is equivalent to $A$
having a central element of trace one.

\begin{lem}\label{ZentralesSpur-1Element} Let $H$ be a Hopf algebra over $R$ with bijective antipode and
let $A$ be a left $H$-module algebra. Then there exists a right integral $t$ in $H$ and a central element $z$ in $A$ with
$S(t)\cdot z=1$ if and only if $A$ has a central element of trace one.
\end{lem}

\begin{proof}
Let $t\in \int_r$ and $z\in Z(A)$ with  $S(t)z=1$. Since the antipode is bijective $S(t)\in \int_l$ and $A$ has a central element of trace one.
On the other hand if $z$ is a central element of trace one, then $t\cdot z=1$ for some left integral $t$.
Thus $S^{-1}(t)$ is a right integral in $H$ such that $S(S^{-1}(t))\cdot z= 1$.
\end{proof}



\subsection{}

The existence of a central element of trace one for group actions implies the separability of the skew group ring.
More generally we have the following theorem
that has also been proved in  \cite[Theorem 1.11]{CohenFishman92}
(see also  \cite{CohenFishman94}) for crossed products but under the additional
hypothesis of $H$ being a Frobenius $R$-algebra. As we will see, we do not need any additional assumptions.

\begin{thm}\label{CocommutativeTraceOneSeparable} Let $H$ be a cocommutative Hopf algebra over $R$ and $A$ a left $H$-module algebra.
If $A$ contains a central element of trace one then $\AH$ is a separable extension of $A$.
\end{thm}

\begin{proof}
Let $z$ be a central element of trace one with $t\cdot z=1$ for some $t\in \int_l$.
Note that since $H$ is cocommutative $S^2=id$.
Hence $t':=S(t)=S^{-1}(t)\in \int_r$ and $S(t')\cdot z = 1$.
Let $\Omega$ be the element in $\AH {\otimes_A} \AH$ from the proof of \ref{SeparabeleSmashProdukte}
with respect to $t'$ and $z$. Then
$$\Omega = \sum_{(t')} (1\# S(t'_1)) \otimes (z\# t'_2) = \sum_{(t)} (1\# t_2) \otimes (z\# S(t_1))$$
is $\AH$-centralising in  $\AH {\otimes_A} \AH$ as it has been
shown in the proof of \ref{SeparabeleSmashProdukte}. Moreover
$$ \mu(\Omega) = \omega =  \sum_{(t)} (t_2\cdot z)\# t_3S(t_1) = \sum_{(t)} (t_1\cdot z) \# t_2S(t_3) = (t\cdot z)\# 1 = 1\# 1$$
where $\mu: \AH{\otimes_A}\AH \rightarrow \AH$ denotes the multiplication of $\AH$. Hence $\Omega$ is a separability idempotent, i.e.
$\AH$ is separable over $A$.
\end{proof}

\subsection{}

The following Corollary shows that for group actions or more generally for actions of cocommutative Hopf algebras the separability of the smash
product and the projectivity of the module algebra coincide.

\begin{cor} Let $H$ be a cocommutative Hopf algebra over $R$ and let $A$ be a commutative left $H$-module algebra
that is projective as $R$-module. The following statements are equivalent:
\begin{enumerate}
 \item[(a)] $\AH$ is a  separable extension of $A$.
 \item[(b)] $A$ is a  projective left $\AH$-module.
 \item[(c)] $A$ contains an element of trace one.
\end{enumerate}
\end{cor}

\begin{proof}
$(a)\Rightarrow (b)$ Since separable extensions are semisimple extensions and since
the $\AH$-linear projection $\AH\rightarrow A$ splits in $A$-Mod it also splits in $\AH$-Mod.
Hence $A$ is projective as left $\AH$-module.\\
$(b)\Leftrightarrow (c)$ Follow from \ref{Spur1Element2}.\\
$(c)\Rightarrow (a)$ Follows from \ref{CocommutativeTraceOneSeparable}.
\end{proof}

\subsection{}\label{BialgebrenNichtSeparabel} $(b)\Rightarrow (a)$ does not hold in general for
bialgebras. For instance let $S$ be the multiplicative monoid
$(\ZZ_3,\cdot, 1)$ and consider the 3-dimensional commutative
cocommutative $\ZZ_2$-bialgebra $H:=\ZZ_2[S]$. One checks easily
that the counit $\varepsilon:H\rightarrow \ZZ_2$ splits as ring
homomorphism and yields a decomposition $H \simeq \ZZ_2 \times
\ZZ_2[C_2]$ as $\ZZ_2$-algebras where $C_2$ denotes the group of
two elements. Note that $\ZZ_2[C_2]$ is not semisimple as
$\ZZ_2$-algebra and hence $H$ is not semisimple and therefore also
not a separable extension of $\ZZ_2$. On the other hand the
trivial $H$-module algebra $\ZZ_2$ is projective as $H$-module.

Hence there exists a cocommutative finite dimensional bialgebra
over a field and a commutative module algebra $A$  such that $A$
is projective as $\AH$-module, but $\AH$ is not a separable
extension of $A$.

%

\subsection{}
Analogous to the characterisation of separable Hopf algebras we can characterise separable smash products.

\begin{thm}\label{SeparabeleSmashProdukteInvertibleEpsilon}
Let $H$ be a Hopf algebra over $R$ and $A$ a left $H$-module algebra. If $H$ contains a left or right integral $t$ such that
$\varepsilon(t)$ is invertible in $A$, then $\AH$ is a separable extension of $A$.
\end{thm}

\begin{proof} Let $t$ be a right integral such that $\varepsilon(t)$ is invertible in $A$ and denote by
$z=\varepsilon(t)^{-1} \in A$ the inverse of $\varepsilon(t)$.
Since $\varepsilon(t)\in R1_A \subseteq Z(A)$, $z$ is a left and right inverse of $\varepsilon(t)$.
Moreover it is easy to see that $z$ is also a central element of $A$.
For all $h\in H$ consider the following:
\[hz-\varepsilon(h)z = z\varepsilon(t) [ hz - \varepsilon(h)z ]
= z[ h(\varepsilon(t)z) - \varepsilon(h)(\varepsilon(t)z) ]
= z[ h1_A - \varepsilon(h)1_A]= 0,
\]
Hence $z$ is an $H$-invariant element. Thus $S(t)\cdot z =
\varepsilon(t)\cdot z = 1$ holds. Let $\Omega$ the element from
the proof of \ref{SeparabeleSmashProdukte} with respect to $t$ and
$z$. Then $\Omega$ is $\AH$-centralising. Moreover
$$\mu(\Omega)=\omega = \sum_{(t)} (S(t_2)\cdot z)\# S(t_1)t_3 =
\sum_{(t)} z\# \varepsilon(t_2)S(t_1)t_3 = z\# \varepsilon(t) =
1\#1$$ shows that $\Omega$ is a separability idempotent for $\AH$,
i.e. $\AH$ is separable over $A$.
\end{proof}

\subsection{}
One should compare the hypothesis of the last theorem with the
(frequent) assumption of the existence of $|G|^{-1}$ in $A$ in the theory of (finite) group actions. Since $t=\sum_{g\in G} g$ is an integral
of the group ring $\ZZ[G]$ and $\varepsilon(t)=|G|$ our last theorem showed that the invertibility of $|G|$ in $A$ implies
the separability of $A*G = A\# \ZZ[G]$ over $A$.

Also the presumable weaker condition that $A$ is $|G|$-torsionfree allows to localise $A$ by the powers of $|G|$ and obtain an
algebra $A'$ where $|G|$ is invertible. Thus one can go over from $A$ to a localisation $A'$ of $A$ such that $A'*G$ is separable over $A'$.
The same argument is indeed applicable to a $H$-module algebra $A$.

\subsection{}
The condition of Theorem \ref{SeparabeleSmashProdukteInvertibleEpsilon} is obviously fulfilled if $\varepsilon(t)$ is invertible in $R$, i.e. if
$H$ is separable over $R$. Hence we get the following extension of characterisation of separable Hopf algebras.

\begin{cor} The following statements are equivalent for a Hopf algebra $H$ over $R$.
\begin{enumerate}
    \item[(a)] $H$ is a separable $R$-algebra.
    \item[(b)] $\AH$ is a separable extension of $A$ for any left $H$-module algebra $A$.
    \item[(c)] Every left $H$-module algebra has an element of trace one.
    \item[(d)] Every left $H$-module algebra $A$ is a projective left $\AH$-module.
\end{enumerate}
\end{cor}

\bibliographystyle{amsalpha}
\providecommand{\bysame}{\leavevmode\hboxto3em{\hrulefill}\thinspace}
\providecommand{\MR}{\relax\ifhmode\unskip\space\fi MR }
\providecommand{\MRhref}[2]{%
  \href{http://www.ams.org/mathscinet-getitem?mr=#1}{#2}
} \providecommand{\href}[2]{#2}

\end{document}